\documentclass[12pt]{amsart}

\usepackage{amsmath,amssymb,amscd}
\usepackage{epsfig}
\usepackage{latexsym}

\theoremstyle{remark}
\newtheorem*{remark}{\bf Remark}
\theoremstyle{plain}
\newtheorem{theorem}{\bf Theorem}
\newtheorem*{theorem*}{\bf Theorem}
\newtheorem{lemma}{\bf Lemma}
\newtheorem{definition}{\bf Definition}
\newtheorem{proposition}{\bf Proposition}
\newtheorem{corollary}{\bf Corollary}

\def\Q{{\mathbb Q}}
\def\Z{{\mathbb Z}}
\def\R{{\mathbb R}}
\def\C{{\mathbb C}}

\def\T{{\mathbb T}}

\def\e{{\varepsilon}}

\def\trans{{\rm Trans}}
\def \eqdef{:=}
\def \cal{\mathcal}
\def\eps{\varepsilon}
\def \Id{\rm Id}

\def\ds{\displaystyle}

\title{On the widths of the Arnol'd Tongues}
\author{Kuntal Banerjee}
\begin{document}

\maketitle

\begin{abstract}
Let $F: \R \to \R$ be a real analytic increasing diffeomorphism with
$F-\Id$ being 1 periodic. Consider the translated family of maps
$(F_t :\R \to \R)_{t\in \R}$ defined as $F_t(x)=F(x)+t$.  Let
${\trans}(F_t)$ be the translation number of $F_t$ defined by:
\[{\trans}(F_t) \eqdef \lim_{n\to +\infty}\frac{F_t^{\circ n}-{\rm
Id}}{n}.\] Assume there is a Herman ring of modulus $2\tau$
associated to $F$ and let $p_n/q_n$ is the $n$-th convergent of
$\trans(F)$. Denoting $\ell_{\theta}$ as the length of the interval
$\{t\in \R~|~\trans(F_t)=\theta\}$, we prove that the sequence
$(\ell_{p_n/q_n})$ decreases exponentially fast with respect to
$q_n$. More precisely
\[\limsup_{n \to \infty} \frac{1}{q_n} \log {\ell_{p_n/q_n}} \le
-2\pi \tau .\]
\end{abstract}

\bigskip

\section{Introduction}

In the whole article, $F:\R\to \R$ is an increasing analytic
diffeomorphism such that $F-\Id$ is 1 periodic. We identify the
circle $\T \eqdef \R/\Z$ and $\mathbb S^1$ via $[x]\simeq e^{2i\pi
x}$, where $[x]$ denotes the class of real number $x$ modulo 1. The
map $F$ induces an orientation preserving analytic circle
diffeomorphism $f:\T \to \T$ given by $[x]\mapsto [F(x)]$. The
definitions of translation number of $F$ and the rotation number of
$f$ are based upon the following result of Poincar\'e.

\begin{theorem*}[Poincar\'e]
The sequence of maps $\ds \frac{F^{\circ n}-{\rm Id}}{n}$ converges
uniformly on $\R$ to a constant.
\end{theorem*}

\begin{definition}
The {\em translation number} of $F$ is defined as \[{\trans}(F)
\eqdef \lim_{n\to +\infty}\frac{F^{\circ n}-{\rm Id}}{n}.\] The {\em
rotation number} of $f$ is the quantity ${\rho}(f) \eqdef
{\trans}(F)\mod 1$.
\end{definition}

Let us now consider the translated family of maps $F_t:\R \to \R$
defined by $F_t(x)=F(x)+t$ for all $x\in \R$ and for every $t \in
\R$. This family induces a family $f_t$ of analytic circle maps. The
function
\[\cal H: t \mapsto \trans(F_t)\] is continuous and non
decreasing. The preimage of an irrational translation number is just
a point under this map. And the preimage of a rational translation
number is a closed interval, generally not reduced to a point. Thus
it is interesting to study the lengths of these intervals
\[I(\theta)\eqdef \{t\in \R~|~\trans(F_t)=\theta\}.\] Let's denote the length of the interval
$I(\theta)$ as $\ell_{\theta}$. The graph of $\cal H$ is like a
devil's staircase generally. The function $\theta \mapsto
\ell_{\theta}$ is highly discontinuous and our aim is to estimate
these lengths under certain conditions.

\begin{figure}
\centerline{\rotatebox{90}{\scalebox{.5}{\includegraphics{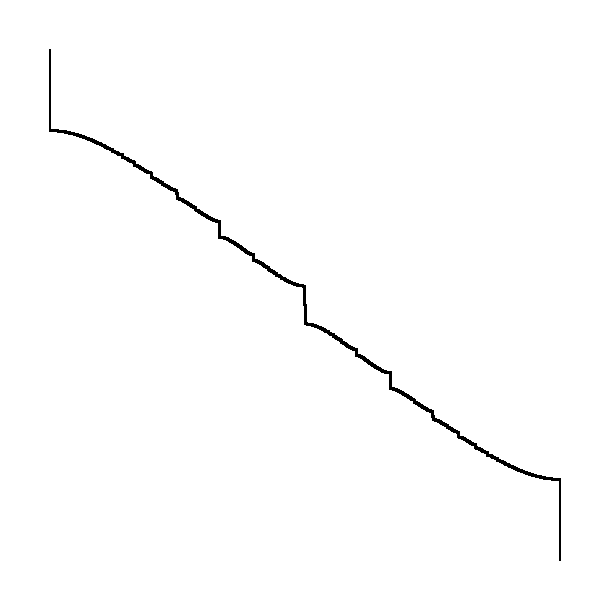}}}}
\caption{The graph of $\cal H: t\mapsto \trans(F_t)$ for $0\le t\le 1$,
where $F_t(x)=x+t+\ds\frac{1}{4\pi}\sin(2\pi x)$.}
\end{figure}

These lengths could be connected with the widths of the Arnol'd
tongues in the following way. Let's define a 2-parameter family of
maps as follows.
\[F_{t,a}(x) = x + t + a\sin(2\pi x).\]
This gives a family of increasing diffeomorphisms of $\R$ when $t\in
\R$ and  $a\in [0,1/2\pi)$. This family is often called the {\em
Arnol'd family} or the {\em standard family} after Arnol'd
\cite{Arnold}.

\begin{definition}
The {\em Arnol'd tongue} $\cal T_\theta$ of translation number
$\theta$ is defined as the following set
\[\cal T_{\theta} \eqdef \{(t,a)\in \R \times [0,1/2\pi)~|~\trans(F_{t,a}) = \theta\}.\]
\end{definition}

If we fix $a\in [0,1/2\pi)$ and set $F_t=F_{t,a}$, then the length
$\ell_{p/q}$ is the width of the Arnol'd tongue $\cal T_{p/q}$
sliced at the height $a$.

\begin{figure}
\centering
\includegraphics[width=5in,height=2.5in]{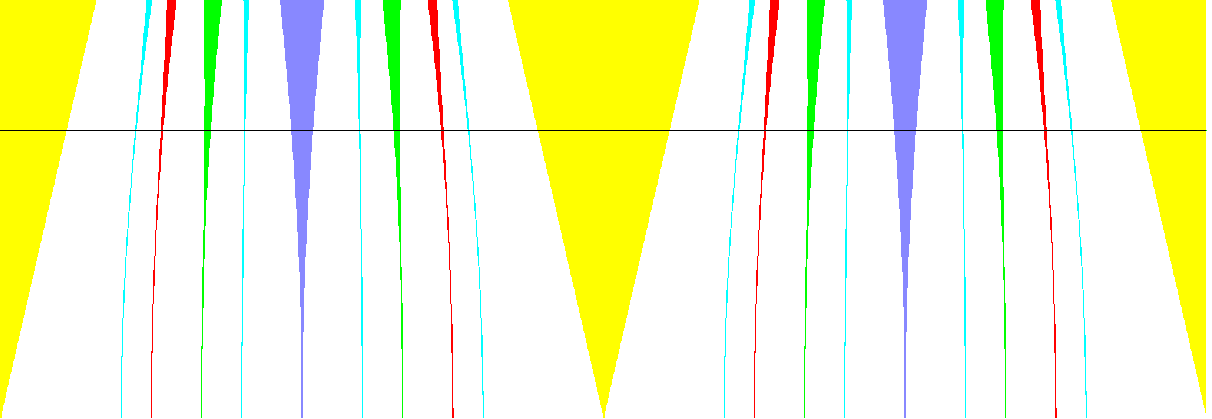}
\caption{Arnol'd tongues of the standard family sliced at a fixed
height.}
\end{figure}

{\bf Acknowledgments :} This work is a part of my thesis, which is
funded by CODY and Marie-Curie Research Training Networks. I would
like to thank my adviser Xavier Buff for giving me this problem and
his guidance. I would also thank Arnaud Ch\'eritat for his picture
of Herman strip.

\section{Preliminaries}

Before we proceed further we would recall some basic facts about
translation and rotation numbers. Every time we write $p/q$ for a
rational number, we implicitely assume that $p$ and $q$ are coprime.

\begin{theorem*}[Poincar\'e]
If $\rho(f)\in \Q/\Z$ then $f:\T\to \T$ has a periodic point. More
precisely, if $\trans(F) =p/q\in \Q$ then there is a point $a\in \R$
such that $F^{\circ q}(a)=a+p$.
\end{theorem*}

Note that $G\eqdef F^{\circ q}-{\rm Id}-p$ vanishes on the whole
$F$-orbit of $a$, in particular on the $q$-set $\{a, F(a), \ldots,
F^{\circ (q-1)}(a)\}$ whose image in $\T$ is a cycle of $f$. We
shall say that such a cycle has rotation number $p/q$. The
derivative of $G$ is constant along the orbit of $a$ under iteration
of $F$. As $G$ is analytic, either it has a double root, or it
vanishes at least once with positive derivative and once with
negative derivative. This shows that counting multiplicities, $f$
has at least $2$ cycles with rotation number $p/q$.

\begin{theorem*}[Poincar\'e]
If $\rho(f) = \alpha \in \R\setminus\Q$, then $f:\T \to \T$ is
semi-conjugate to the rotation $\T \ni [x]\mapsto [x+\alpha] \in
\T$.
\end{theorem*}

In fact, the semiconjugacy may be obtained as follows; the sequence
of maps
\[\frac{1}{N}\sum_{k=0}^{N-1} \left(F^{\circ k}-k\alpha\right)\]
converges, as $N\to +\infty$, to a non-decreasing continuous
surjective map $\Phi_{F}:\R\to \R$, which satisfies
\[\Phi_{F}(x+1) = \Phi_{F}(x)+1\quad\text{and}\quad \Phi_{F}\circ F(x) =
\Phi_{F}(x) + \alpha.\]

The following result of Denjoy implies that when $F$ is an analytic
diffeomorphism, then the semiconjugacy is in fact an actual
conjugacy. In other words $\Phi_{F}:\R\to \R$ is an increasing
homeomorphism.

\begin{theorem*}[Denjoy]
If $\rho(f)=\alpha \in \R\setminus\Q$ and if $f$ is a ${\cal C}^2$ diffeomorphism,
then $f:\T\to \T$ is conjugate to the rotation of angle $\alpha$.
\end{theorem*}

\section{Herman Ring}

From now on, we assume that $\trans(F)\in \R\setminus\Q$ and so, that
$\Phi_F:\R\to \R$ is a homeomorphism. We will now be interested in
the regularity of $\Phi_F$. It is known that when $\trans(F)$
satisfies an appropriate arithmetic condition, then the conjugacy
$\Phi_F$ is itself an analytic diffeomorphism.

The first result obtained in this direction is a result of Herman.
Recall that $\alpha \in \R\setminus\Q$ is a {\em Diophantine number}
if there are constants $C>0$ and $\tau\geq 2$ such that
\[|\alpha-p/q|\geq C/q^\tau\]
for all rational number $p/q$.

\begin{theorem*}[Herman]
If $\trans(F)$ is a Diophantine number, then $\Phi_{F}:\R\to \R$ is
an analytic diffeomorphism.
\end{theorem*}

For a proof, one can look at \cite{Herbig}.

The optimal arithmetic condition which guaranties that the conjugacy
is an analytic diffeomorphism, has been obtained by Yoccoz
\cite{yoccoz} but is too complicated to be recalled here.

Now, we can introduce the definition of the Herman ring.

\begin{definition}
Assume that $\trans(F)\in \R\setminus\Q$ and $\Phi_{F}:\R\to \R$ is
analytic. Let $\tau$ be the largest number such that $\Psi_{F}
\eqdef \Phi_{F}^{-1}$ extends univalently to $S(\tau) \eqdef
\{z\in\C~|~-\tau<{\rm Im}(z) < \tau\}$. The map $F$ extends
analytically to $\cal{HS}(F)\eqdef \Psi_{F}\bigl(S(\tau)\bigr).$ We
call $\cal{HS}(F)$ the {\em Herman strip} of $F$. The image of
$\cal{HS}(F)$ in $\C/\Z$ is called a {\em Herman ring} of $f$
associated to $F$. The modulus of the Herman ring is $2\tau$.
\end{definition}

\begin{figure}
\centering
\includegraphics[width=4.5in,height=4.5in]{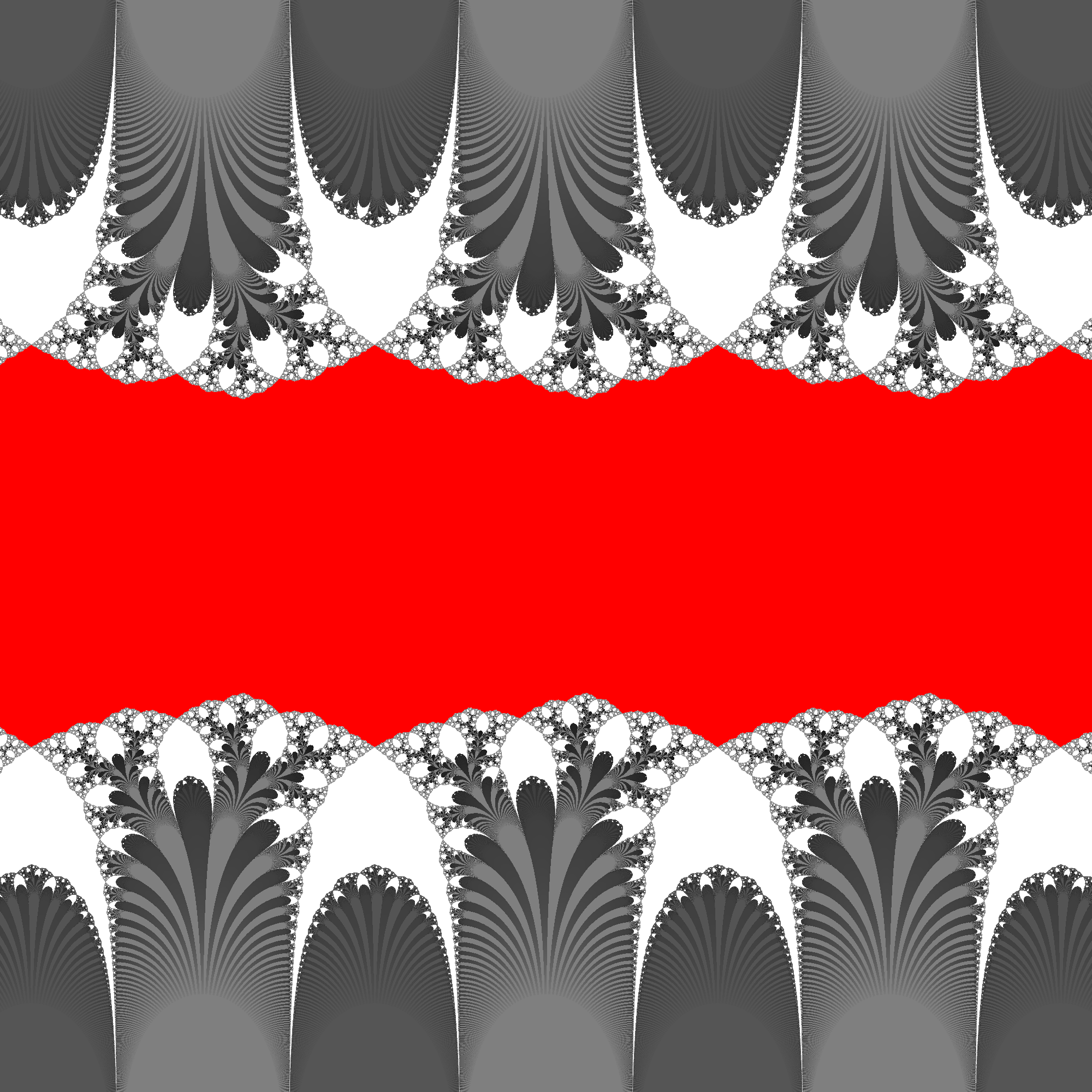}
\caption{A Herman strip in the family
$F_t(x)=x+t+\ds\frac{1}{4\pi}\sin(2\pi x)$. The translation number
is the golden mean $(\sqrt 5-1)/2$.}
\end{figure}

From now on, we assume that $\theta\eqdef \trans(F)$ is irrational
and that $\Phi_F$ is an analytic diffeormorphism, i.e., $F$
has a Herman Strip. Then, we would study the length $\ell_{p_n/q_n}$
where $p_n/q_n$ is the $n$-th convergent of the continued fraction
expansion of $\theta$.

\section{Main result and its comparison to earlier works}

\begin{theorem}
Suppose that $F:\R\to \R$ is an increasing $\R$-analytic
diffeomorphism. For any $t\in \R$, define the translated family of
maps $F_t (x)= F(x)+t$ for $x\in \R$. Assume that
\begin{itemize}
\item $\trans(F) = \theta\in \R\setminus\Q$ and

\item There is a Herman ring associated to $F$ with modulus $2\tau$.
\end{itemize}
Let $p_n/q_n$ be the $n$-th continued fraction convergent of
$\theta$. Then, we have the following inequality
\[\limsup_{n\to \infty} \frac{1}{q_n} \log \ell_{p_n/q_n} \le -2\pi\tau.\]
\end{theorem}

In our set up, when we approach a map which has a Herman strip then
the corresponding lengths $\ell_{p_n/q_n}$ decreases exponentially
with respect to $q_n$. In particular, when we take a horizontal
slice of the Arnol'd tongues, and when we approach a parameter with
a Herman ring associated to it, the width of the tongue
$\ell_{p_n/q_n}$ decreases exponentially with respect to $q_n$.

Herman studied the function $\cal H : t \mapsto \trans(F_t)$ in his
paper \cite{Hershort}. From his works one can have an estimate on
the behaviour of these lengths.

\begin{theorem*}[Herman]
If  $\Phi_{F}$ is ${\cal C}^1$, then the function $\cal H: t\mapsto
\trans(F_t)$ has a non-zero derivative at $t=0$.
\end{theorem*}

\begin{corollary}
If  $\Phi_{F}$ is ${\cal C}^1$, then
\[\ell_{p_n/q_n} = o\left(\theta-\frac{p_n}{q_n}\right) = \cal{O}\left(\frac{1}{q_n^2}\right).\]
\end{corollary}

\begin{proof}
Suppose \[I(p_n/q_n)=[t_n^-,t_n^+].\] The Herman's theorem implies that
\begin{align*}
\frac{p_n/q_n-\theta}{t_n^\pm-0} & \to \cal H'(0)\ \ (\ne 0)\\
\frac{t_n^\pm}{p_n/q_n-\theta} & \to \frac{1}{\cal H'(0)}\\
\Rightarrow t_n^\pm &= \frac{1}{\cal H'(0)}\left(\frac{p_n}{q_n}-\theta\right)+o\left(\frac{p_n}{q_n}-\theta\right).\\
\therefore \ell_{p_n/q_n}&=t_n^+-t_n^-
=o\left(\frac{p_n}{q_n}-\theta\right)=\cal{O}\left(\frac{1}{q_n^2}\right).\qedhere
\end{align*}
\end{proof}

Thus according to Herman $\ell_{p_n/q_n}$ decreases at least as fast
as $1/{q_n^2}$. Our theorem states that the decay is much faster, it
is an exponential decay. More precisely
\[\ell_{p_n/q_n}\le e^{-2\pi \tau q_n + o(q_n)}.\]
This estimate is better than previous works and it involves the
modulus of the Herman ring in the inequality. Since we can have a
Herman ring whenever the rotation number is of bounded type and
satisfies {\em Herman condition} \cite{yoccoz}, this estimate is
valid for a big subset of irrational rotation numbers.

\section {Proof of the main result}
We shall start with estimating the length $\ell_{p/q}$ under some
conditions.

\begin{lemma}
Assume that there are $\eps_0>0$ and $v_0>0$ such that for all $t\in
I(p/q)$ and for all $ x\in \R$,
\begin{itemize}
\item $m_t\le F_t^{\circ q}(x)-x-p \le M_t$, with $M_t-m_t\le \eps_0$
  and
\item $\ds \frac{\partial F_t^{\circ q}(x)}{\partial t}\ge v_0$.
\end{itemize}
Then,
\[\ell_{p/q}\le \frac{\eps_0}{v_0}.\]
\end{lemma}

\begin{proof} Let $I(p/q)=[t_{p/q}^-,t_{p/q}^+]$. As we are in an increasing family, we have $M_{t_{p/q}^-}=0$ and $m_{t_{p/q}^+}=0$.
The assumption that $\displaystyle\frac{\partial F_t^{\circ
q}(x)}{\partial t} \ge v_0$ implies the following.
\begin{align*}
F_{t_{p/q}^+}^{\circ q}(x)-x-p &\ge F_{t_{p/q}^-}^{\circ q}(x)-x-p+v_0(t_{p/q}^+-t_{p/q}^-)\\\Rightarrow m_{t_{p/q}^+} &\ge m_{t_{p/q}^-} + v_0(t_{p/q}^+-t_{p/q}^-) \\
\Rightarrow M_{t_{p/q}^-}-m_{t_{p/q}^-} &\ge v_0(t_{p/q}^+-t_{p/q}^-) \;\ (\because M_{t_{p/q}^-}=m_{t_{p/q}^+}=0)\\
\Rightarrow \frac{\e_0}{v_0} &\ge (t_{p/q}^+-t_{p/q}^-).
\end{align*}
Hence $\displaystyle \bigl|I(p/q)\bigr| \le \frac{\e_0}{v_0}$.
\end{proof}

\begin{remark}
It is easy to see by induction on $k$ that in our family, for all
$t\in \R$, for all $x\in \R$ and for all $k\geq 1$,
\[\frac {\partial F_t^{\circ k}(x)}{\partial t}\geq 1.\]
\end{remark}

This lemma gives an estimate of $I(p/q)$ assuming that $|M_t-m_t|$
is bounded and $\ds \frac{\partial F_t^{\circ q}(x)}{\partial t}\ge
v_0$. According to the remark $v_0$ can be taken as 1 in our family.
Now we are interested to find a bound of $|M_t-m_t|$.

Choose a sequence $\bigl(t_n\in I(p_n/q_n)\bigr)_{n\geq 1}$. Define
 \[G_n(x) \eqdef F_{t_n}^{\circ q_n}(x)-x-p_n.\]
 $G_n$ vanishes along at least two sets of $q_n$ points corresponding to two cycles of period $q_n$ for $f_{t_n}$.

It is enough to show that for all $\tau'<\tau$, we have
\[\sup_{x\in \R} \bigl|G_n(x)\bigr|  =\cal{O}\bigl
(e^{-2\pi \tau' q_n}\bigr).\]

Choose $\tau'<\tau$ and set $S' \eqdef
\Phi_{F}^{-1}\bigl(S(\tau')\bigr)\subset \cal{HS}(F)$.

\begin{proposition}
If $n$ is large enough, then for all $k\leq q_n$, $F_{t_n}^{\circ
k}$ is defined in $S'$ with values in $\cal{HS}(F)$. The sequence of
maps $(G_n)$ converges uniformly to $0$ on $S'$.
\end{proposition}

To prove this proposition let's prove parts of it in the following
two lemmas. Choose $\tau''$ with $\tau'<\tau''<\tau$. Set $S''\eqdef
\Phi_{F}^{-1} \bigl(S(\tau'')\bigr)\subset \cal{HS}(F)$. Note that
$S'\subset S''\subset \cal{HS}(F)$ and $S(\tau')\subset
S(\tau'')\subset S(\tau)$.

\begin{lemma}
$\displaystyle \sup_{S''}|F^{\circ q_n}(z)-z-p_n| \rightarrow 0$ as
$n \rightarrow \infty$.
\end{lemma}

\begin{proof}
Suppose $w=\Phi_{F}(z)$, then $\Phi_{F}\circ F^{\circ q_n}(z)=w+q_n
\theta$. By choice $(\Phi_{F}^{-1})'$ is bounded on
$\overline{S(\tau'')}$. Assuming $C_1=\displaystyle \max_{w\in
\overline{S(\tau'')}}|(\Phi_{F}^{-1})^{\prime}(w)|$, we have
\begin{align*}
|F^{\circ q_n}(z)-z-p_n| & =|\Phi_{F}^{-1}(w+q_n \theta -p_n)-\Phi_{F}^{-1}(w)|\\
 & \le C_1|w+q_n \theta - p_n-w| \;\ (\mbox{by Mean Value Inequality})\\
 & = C_1|q_n \theta - p_n| \\
 & \le \frac{C_1}{q_n}.
\end{align*}
The last inequality uses the fact that $p_n/q_n$ is the $n$-th
continued fraction convergent of $\theta$, thus $|\theta-p_n/q_n|\le
1/q_n^2$. Hence
\[\sup_{S''} |F^{\circ q_n}(z)-z-p_n|\underset{n\to
+\infty}\longrightarrow 0.\qedhere\]
\end{proof}

Before we start next lemma define $H_n: w \mapsto \Phi_{F}\circ
F_{t_n} \circ \Phi_{F}^{-1}$ on $S(\tau')$.

\begin{lemma}
For $n$ large enough and for all $k \le q_n$, $H_n^{\circ k}$ is
defined on $S(\tau')$ and $H_n^{\circ k}S(\tau') \subset S(\tau)$.
Moreover
\[\sup_{w \in S(\tau')} |H_n^{\circ q_n}(w)-w-q_n \theta| \rightarrow 0.\]
\end{lemma}

\begin{proof}
We have to show that for large
$n$ and for all $k \le q_n$, $H_n^{\circ k}$ is defined on
$S(\tau')$ and the image of $H_n^{\circ k}$ is inside $S(\tau)$. For
$z\in S'$ and $k\le q_n$ we have
\begin{align*}
|\Phi_{F}\circ F_{t_n}^{\circ k}(z)-\Phi_{F}\circ F^{\circ k}(z)| &=  |\Phi_{F}\circ F_{t_n}^{\circ k}(z)-\Phi_{F}(z)-k \theta|  \\
&=  \left| \sum_{j=1}^k\Phi_{F}\circ F_{t_n}^{\circ j}(z)-\Phi_{F}\circ F_{t_n}^{\circ j-1}(z)-\theta\right|\\
& \le \sum_{j=1}^k\bigl|\Phi_{F}\circ F_{t_n}^{\circ j}(z)-\Phi_{F}\circ F_{t_n}^{\circ j-1}(z)-\theta\bigr|\\
& = \sum_{j=1}^k\bigl|\Phi_{F}\circ F_{t_n}^{\circ
j}(z)-\Phi_{F}\circ F\circ F_{t_n}^{\circ j-1}(z)\bigr|
\end{align*}

Assume that for any $j\le q_n$ the point $F_{t_n}^{\circ j}(z)$ is
inside $S''$ for $z\in S'$. Also set $C_2=\ds\max_{z\in
\overline{S''}}|\Phi_{F}'(z)|$. As $|F_{t_n}(z)-F(z)|= |t_n|$ for
any $z$ we see that
\[ |\Phi_{F}\circ F_{t_n}^{\circ k}(z)-\Phi_{F}\circ F^{\circ k}(z)|  \le C_2 k|t_n| \;\ (\text {by Mean Value Inequality}).\]

Since $\cal H: t\mapsto \trans (F_t)$ has a non zero derivative at
$0$, we have
\[|\theta-\frac{p_n}{q_n}| \ge |\cal H'(0)/2||t_n|.\]
Choosing $C_3=\ds\frac{2}{|\cal H'(0)|}$ we see that

\begin{align*}
|F_{t_n}^{\circ k}(z)-F^{\circ k}(z)| &=
\bigl|\Phi_{F}^{-1}\bigl(\Phi_{F}\circ F_{t_n}^{\circ
k}(z)\bigr)-\Phi_{F}^{-1}
\bigl(\Phi_{F}\circ F^{\circ k}(z)\bigr)\bigr|\\
& \le C_1 \bigl|\Phi_{F}\circ F_{t_n}^{\circ k}(z)-\Phi_{F}\circ F^{\circ k}(z)\bigr| \\
&   \;\ (\mbox{by Mean Value Inequality})\\
& \le C_1C_2 k|t_n|  \\
& \le C_1C_2C_3\frac {k}{q_n^2}.
\end{align*}
In the last inequality we use the fact that $(p_n/q_n)$ are the
convergents to $\theta$. Thus if $n$ is large enough the point
$F_{t_n}^{\circ k}(z)$ is inside $S''$ when $z$ is taken in $S'$ for
all $k\le q_n$. This means our assumption above is true for large
$n$. Consequently for large $n$ and for all $k \le q_n$, the
conjugate $H_n^{\circ k}$ is defined on $S(\tau')$ and the image of
$H_n^{\circ k}$ is inside $S(\tau'')\subset S(\tau)$. Taking
$w=\Phi_{F}(z)$ and from the previous calculations we see that
\[|H_n^{\circ q_n}(w)-w-q_n \theta|=\bigr|\Phi_{F}\circ F_{t_n}^{\circ q_n}(z)-\Phi_{F}\circ F^{\circ q_n}(z)\bigr|
 \le \displaystyle \frac{C_2C_3}{q_n} .\]
 And hence $\ds \sup_{w \in S(\tau')} |H_n^{\circ q_n}(w)-w-q_n \theta| \rightarrow 0$ as $n \to \infty$.
\end{proof}

\begin{proof}[Proof of the Proposition 1:]
In Lemma 3 we have already seen that for large $n$ and for all
$k\leq q_n$, $F_{t_n}^{\circ k}$ is defined in $S'$ with values in
$\cal{HS}(F)$. And for $z\in S'$,
\begin{align*}
\bigl|G_n(z)\bigr|=\bigl|F_{t_n}^{\circ q_n}(z)-z-p_n\bigr| & = \bigl|F_{t_n}^{\circ q_n}(z)-F^{\circ q_n}(z)+F^{\circ q_n}(z)-z-p_n\bigr|\\
& \le \bigl|F_{t_n}^{\circ q_n}(z)-F^{\circ q_n}(z)\bigr| +
\bigl|F^{\circ q_n}(z)-z-p_n\bigr|.
\end{align*}
From Lemma 2 we know that $\displaystyle \sup_{S''} \bigl|F^{\circ
q_n}(z)-z-p_n\bigr| \rightarrow 0$ as $n\to \infty$. And
\begin{align*}
|F_{t_n}^{\circ q_n}(z)-F^{\circ q_n}(z)|&=|\Phi_{F}^{-1}\bigl(\Phi_{F}F_{t_n}^{\circ q_n}(z)\bigr)-\Phi_{F}^{-1}\bigl(\Phi_{F}F^{\circ q_n}(z)\bigr)|\\
& \le C_1 |\Phi_{F}F_{t_n}^{\circ q_n}(z)-\Phi_{F}F^{\circ q_n}(z)| \\
& \;\ (\mbox{by Mean Value Inequality})\\
& \le \frac{C_1C_2C_3}{q_n} \;\ (\mbox{by Lemma 3}).
\end{align*}
This completes the proposition.
\end{proof}

We set \[\Phi_n \eqdef
\frac{1}{q_n}\sum_{k=0}^{q_n-1}
( F^{\circ k} - k\theta)\quad\text{and}\quad
\widehat \Phi_n \eqdef \frac{1}{q_n}\sum_{k=0}^{q_n-1}
 \left( F_{t_n}^{\circ k} - k\frac{p_n}{q_n}\right).\]
We know that
\[\Phi_n\underset{n\to +\infty}\longrightarrow \Phi_F.\]
We shall see that this is also true for the sequence $\widehat
\Phi_n$.

\begin{lemma}
If $n$ is large enough, the domain of $\widehat \Phi_n$ contains
$S'$. As $n\to \infty$, the sequence $\widehat \Phi_n$ converges to
the linearizing map $\Phi_{F}$ and thus $\widehat \Phi_n$ has a
univalent inverse $\widehat \Psi_n : S(\tau')\to \cal{HS}(F)$ for
large $n$.
\end{lemma}

\begin{proof}
Taking $z\in S'$ we note that for any $k\le q_n$ the point
$F_{t_n}^{\circ k}(z)$ is inside $S''$ for large $n$ and,
\begin{align*}
\left| \widehat \Phi_{n}(z)-\Phi_{n}(z)\right| & = \frac{1}{q_n}\left|\sum_{k=0}^{q_n-1} \bigl(F_{t_n}^{\circ k}(z)- F^{\circ k}(z)\bigr)-k\left(\theta-\frac{p_n}{q_n}\right)\right|\\
& \le \frac{1}{q_n^3}\sum_{k=0}^{q_n-1} k(C_1C_2C_3-1)\;\ (\mbox {by Lemma 3})\\
& = \frac{1}{q_n^3}\frac{q_n(q_n-1)}{2}(C_1C_2C_3-1)\\
& = \frac{q_n-1}{2q_n^2}(C_1C_2C_3-1).
\end{align*}
This implies that as  $n\to \infty$, the sequence $\widehat \Phi_n$
converges to the linearizing map $\Phi_{F}$. Moreover $\widehat
\Phi_n : S' \to S(\tau)$ is defined and it is univalent for large
$n$. Thus it has a univalent inverse $\widehat \Psi_n: S(\tau')\to
\cal{HS} (F)$ for large $n$.
\end{proof}

\begin{lemma}
Counting multiplicities, the map $G_n\circ \widehat \Psi_n$ vanishes
at least along two sets of the form $a_n + k p_n/q_n$ and $b_n + k
p_n/q_n$ with $a_n\in \R$, $b_n\in \R$ and $k\in \Z$.
\end{lemma}

\begin{proof}
The map $G_n$ vanishes on two sets of $q_n$ points on $\R/\Z$,
counted with multiplicity, corresponding to two $q_n$-cycles of
$f_{t_n}$. Lets assume that $G_n$ vanishes on the sets
$\{F_{t_n}^{\circ j}(a)+k~|~k\in \Z\}$ and $\{F_{t_n}^{\circ
j}(b)+k~|~k\in \Z\}$ for $j=0,\cdots,q_n-1$ and for some $a,b\in
\R$. We have
\[\widehat \Phi_n(a)= \frac{1}{q_n}\sum_{k=0}^{q_n-1}
 \left( F_{t_n}^{\circ k}(a) - k\frac{p_n}{q_n}\right).\] And
 \begin{align*}
 \widehat \Phi_n\bigl(F_{t_n}(a)\bigr) &= \frac{1}{q_n}\sum_{k=0}^{q_n-1}
 \left( F_{t_n}^{\circ k+1}(a) - k\frac{p_n}{q_n}\right) \\
& =  \frac{1}{q_n}\sum_{k=0}^{q_n-1}
 \left( F_{t_n}^{\circ k+1}(a) - (k+1)\frac{p_n}{q_n}+\frac{p_n}{q_n}\right)\\
& = \widehat \Phi_n(a)+\frac{p_n}{q_n}
\end{align*}
In similar way we can obtain
\begin{align*}
& \widehat \Phi_n \bigl(F_{t_n}^{\circ j+1}(a)\bigr)-\widehat \Phi_n
\bigl(F_{t_n}^{\circ j}(a)\bigr)\\ &= \widehat \Phi_n
\bigl(F_{t_n}^{\circ j+1}(b)\bigr)-\widehat \Phi_n
\bigl(F_{t_n}^{\circ j}(b)\bigr)=\ds\frac{p_n}{q_n}
\end{align*}
for all $0\le j \le q_n-2$. This proves that $G_n \circ \widehat
\Psi_n$ vanishes at least along two sets of the form $a_n + k
p_n/q_n$ and $b_n + k p_n/q_n$ with $a_n\in \R$, $b_n\in \R$ and
$k\in \Z$.

\end{proof}

Now note that the map $G_n$ and $\widehat \Psi_n$ are 1 periodic.
For $n$ large enough, the map $G_n\circ \widehat \Psi_n$ is
holomorphic on $S(\tau')$. And as $n\to \infty$, the sequence
$(G_n\circ \widehat \Psi_n)$ converges uniformly to $0$ on
$S(\tau')$. The required estimate is a consequence of the following
lemma.

\begin{lemma}
Assume $G$ is holomorphic on $S(\tau')$, $1$ periodic and vanishes
on two sets of the form $a + k p/q$ and $b + k p/q$ with $k\in
  \Z$.  Then,
\[\max_{x\in \R}\bigl|G(x)\bigr|\leq \frac{\ds\sup_{z\in
    S(\tau')}\bigl|G(z)\bigr|}{\bigl({\rm sinh}(\pi q \tau')\bigr)^2}
\underset{q\to +\infty}=\cal{O}\left(\sup_{z\in
    S(\tau')}\bigl|G(z)\bigr|\cdot e^{-2\pi q\tau'}\right).\]
\end{lemma}

\begin{proof}
Suppose $G_0(z)=\sin\bigl(\pi q (z-a)\bigr) \sin\bigl(\pi q
(z-b)\bigr)$. Then $G_0$ vanishes exactly on the set of the form
$a+kp/q$ and $b+kp/q$ with $k\in \Z$. Hence the function
$\ds\frac{G(z)}{G_0(z)}$ does not have a pole and thus it is
holomorphic on the strip $S(\tau')$. By Maximum Modulus Principle,
for $z\in S(\tau')$ we have
\[ \left|\frac {G(z)}{G_0(z)}\right| \le  \frac{\ds \sup_{z\in S(\tau')}|G(z)|}{\displaystyle \inf_{z\in \partial S(\tau')}|G_0(z)|}. \]
Since $G$ and $G_0$ are non constant and periodic, these supremum
and infimum values actually occur at the boundary of $S(\tau')$.

Assuming $z=x+i\tau'$ we see that
\begin{align*}
|\sin \bigl(\pi q(z-a)\bigr)| & = \Bigl|\frac{e^{i\pi q(z-a)}-e^{-i\pi q(z-a)}}{2i}\Bigr|\\
& = \Bigl|\frac{e^{i\theta}}{i}\Bigr| \Bigl|\frac{e^{-\pi
q\tau'}-e^{-\pi q\tau'}e^{-2i\theta}}{2}\Bigr|
\;\ (\mbox {where} \;\ \theta=\pi qx-\pi qa)\\
& \ge \frac{e^{\pi q\tau'}}{2}-\frac{e^{-\pi q\tau'}}{2}\\
& = \sinh (\pi q\tau').
\end{align*}
For $q$ large enough, such that $\pi q\tau' > \sqrt{\log 2}$, we see
that
\[|\sin \bigl(\pi q(z-a)\bigr)| \ge \sinh (\pi q\tau')> \frac{e^{\pi q \tau'}}{4}.\]
This proves the lemma.
\end{proof}

\begin{proof}[Proof of the Theorem 1:]
By Lemma 6,
\begin{align*}
\max_{y\in \R}\bigl|G_n \circ \widehat \Psi_n (y)\bigr|& \leq \frac{\ds\sup_{w\in S(\tau')}\bigl|G_n\circ \widehat \Psi_n (w)\bigr|}{\bigl({\rm sinh}(\pi q \tau')\bigr)^2}\\
& \underset{q_n\to +\infty} =\cal{O}\left(\sup_{w\in
    S(\tau')}\bigl|G_n\circ \widehat \Psi_n (z)\bigr|\cdot e^{-2\pi q_n\tau'}\right).
\end{align*}
Since $\widehat \Psi_n$ is invertible and by Lemma 2 and 3, $\ds
\sup_{z \in S'}|G_n(z)| \to 0$ as $n \to \infty$, we see that
\[\sup_{x\in \R}|G_n(x)|=\cal O\bigl(e^{-2\pi \tau' q_n}\bigr)\]
for any $\tau' < \tau$. Therefore the theorem is proved.
\end{proof}

\bigskip

{\em Email Address :} \verb+kuntalb@gmail.com+

Universit\'e de Toulouse (UT3,CNRS,INSA,UT1,UT2), Institut de Math\'ematiques de
Toulouse, 118, route de Narbonne, 31062 Toulouse, France.

\end{document}